\documentclass[11pt, a4paper, oneside]{amsart} 

\usepackage{amsmath, amssymb, amsthm, amsrefs}
\usepackage{mathrsfs}
\usepackage[italicdiff]{physics}
\usepackage{mathtools}

\usepackage{hyperref}
\hypersetup{
    colorlinks=true, 
    linktoc= page,   
    citecolor=blue
}

\allowdisplaybreaks[4]
\numberwithin{equation}{section}

\def\le{\leqslant}
\def\ge{\geqslant}

\def\N{\mathbb{N}}

\theoremstyle{plain}
\newtheorem{theorem}{Theorem}[section]

\newtheorem{lemma}[theorem]{Lemma}

\newtheorem{proposition}[theorem]{Proposition}

\theoremstyle{remark}

\begin{document}
\title[Correction to the article ``GWP and scattering in weighted space for NLS'']{Correction to the article ``Global well-posedness and scattering in weighted space for nonlinear Schr\"{o}dinger equations below the Strauss exponent without gauge-invariance''}
\author[M. Kawamoto]{Masaki Kawamoto}
\address[]{Research Institute for Interdisciplinary Science, Okayama University, Okayama City Okayama, 700-8530, Japan}
\email{kawamoto.masaki@okayama-u.ac.jp}
\author[S. Masaki]{Satoshi Masaki}
\address[]{Department of Mathematics, Hokkaido University, Sapporo Hokkaido, 060-0810, Japan}
\email{masaki@math.sci.hokudai.ac.jp}
\author[H. Miyazaki]{Hayato Miyazaki}
\address[]{Teacher Training Courses, Faculty of Education, Kagawa University, Takamatsu, Kagawa 760-8522, Japan}
\email{miyazaki.hayato@kagawa-u.ac.jp}
\date{}

\maketitle
\begin{abstract}
This article resolves some errors in the paper
``Global well-posedness and scattering in weighted space for nonlinear Schr\"{o}dinger equations below the Strauss exponent without gauge-invariance, Math. Ann. 392, 1051--1097 (2025)''.
The errors are in the proof of contraction of a map associated with our equation in two and three dimensions.
\end{abstract}

\section{Errors}
In this paper we resolve some errors in \cite{KMM25}.
The issue lies in the construction of the contraction mapping used to solve the integral equation \cite[(1.19) in Section 3]{KMM25}.
More precisely, the error occurs in the proof of \cite[(3.12) in Lemma 3.4]{KMM25}.
To obtain
\begin{align*}
	\norm{A_{4,n}(\cdot,v_1)-A_{4,n}(\cdot,v_2)}_{X} \lesssim n^{\frac{7}2}M^{\alpha-1} \norm{v_1-v_2}_{X}
\end{align*}
for $v_1,v_2 \in Z_M$, we need the estimate
\begin{align}
    \norm{ p\left(\frac{\partial F_n}{\partial z} (v_1) -\frac{\partial F_n}{\partial z} (v_2) \right) }_{L^{q_{0}}((0,1];L^2)}
    \lesssim n M^{\alpha -1} \norm{v_1-v_2}_{X},
    \label{iss:1}
\end{align}
where all notations are the same as in \cite{KMM25}.

However, a direct calculation yields
\begin{align*}
    \abs{p\left(\frac{\partial F_n}{\partial z} (v_1) -\frac{\partial F_n}{\partial z} (v_2) \right)}
    \lesssim |v_{1}|^{\alpha-2} |p (v_{1} - v_{2})|
    + \left( |v_{1}|^{\alpha-3} + |v_{2}|^{\alpha-3} \right) |p v_{1}| |v_{1} - v_{2}|,
\end{align*}
which shows that the restriction $\alpha \ge 3$ is required to establish \eqref{iss:1}, unless $F_{n}$ is smooth enough, for instance $F_n(u) = |u|^{2-n}u^{n}$ with $n \in 2 \N$ for $d=3$.

Hence, under the general assumption \cite[Assumption 1.1]{KMM25},
the proof of \cite[Lemma 3.4]{KMM25} is incomplete for $\alpha <3$,
and consequently the main result \cite[Theorem 1.4]{KMM25} is not justified in dimensions $d =2, 3$,
since the admissible range of $\alpha$ is given by
\begin{align*} 
\alpha \in \left( \frac72,5 \right) \quad (d=1), \qquad
\alpha \in (2,3) \quad (d=2), \qquad
\alpha = 2 \quad (d=3).
\end{align*}

The implications of this issue depend on the dimension $d$:
\begin{itemize}
    \item For $d=1$, the proof remains valid since the assumption on $\alpha$ implies $\alpha > 3$.
    \item For $d=3$, where $\alpha=2$, the proof holds for smooth nonlinearities such as $F(u)=u^2$ but fails for non-smooth ones, e.g., $F(u) = |u|^2u/(|u|-au)$ with $a \in (-1,1) \setminus \{0\}$.
    \item For $d=2$, which covers the range $2 < \alpha < 3$, the proof requires modification for all nonlinearities considered in \cite[Assumption 1.1]{KMM25}.
\end{itemize}
The error occurs in the proof of the contraction property in the contraction mapping argument, specifically in the estimate of the difference.
On the other hand, the boundedness part of the argument remains valid.
Thus, it may still be possible to recover the result by working in a smaller function space, for instance $H^1 \cap H^{0,1}$, combined with a compactness argument.

The original purpose of \cite{KMM25} was to introduce a method for treating scattering solutions, applicable to nonlinearities with fractional powers, based on a contraction mapping argument after a suitable transformation of the problem.
While such an alternative approach based on compactness may recover the result, it is desirable to establish it within a direct contraction framework.

In the present paper, we correct the argument so that the main theorem holds for the full class of nonlinearities in \cite[Assumption 1.1]{KMM25}.
As a result, the function space used for the estimates becomes slightly more involved, while the argument itself is streamlined.



\section{Correction}

\subsection{Reformulation of the integral equation and the solution space}

To fix this error, we reformulate the integral equation \cite[(1.19)]{KMM25}.
To formulate \cite[(1.19)]{KMM25}, we employ \cite[(1.8)]{KMM25} to have the identity
\begin{align*}
	i\frac{d}{ds} F_n(v)
	={}&\frac{\partial F_n}{\partial z} (v)
	\left(\frac{p^2}2 v + s^{\alpha_0-2} \tilde{F}(s,v)\right)
	- \frac{\partial F_n}{\partial \overline{z}} (v) \overline{\left(\frac{p^2}2 v + s^{\alpha_0-2} \tilde{F}(s,v) \right)}.
\end{align*}
Instead of the identity, we utilize the following identity:
\begin{align*}
	i\frac{d}{ds} F_n(v)
	={}& \frac{\partial F_n}{\partial z} (v)
	\frac{p^{2}}{2}v
    + \frac{\partial F_n}{\partial z} (v)
	\left( i\partial_{t}v - \frac{p^{2}}{2}v \right)
	- \frac{\partial F_n}{\partial \overline{z}} (v) \overline{\frac{p^{2}}{2}v}
    - \frac{\partial F_n}{\partial \overline{z}} (v) \left(\overline{i\partial_{t}v - \frac{p^{2}}{2}v} \right).
\end{align*}
The identity leads to the integral equation
\begin{align}
    \begin{aligned}
    v(t) ={}& U(t-1) v(1) - i\overline{\lambda_1} I_1(t) \\
    &-i \sum_{n =2}^\infty \overline{ \lambda_n} (A_{1,n}(t) + A_{2,n}(t) + A_{3,n}(t)  + A_{4,n}(t) + \widetilde{A_{5,n}}(t)),
    \end{aligned}
    \label{E:IE}
\end{align}
where $v(1)=e^{i \frac{|x|^2}2} \overline{u_0} \in H^1$,
$I_{1}(t)$ and $A_{j,n}(t)$ ($j=1$, $2$, $3$, $4$) are defined as in \cite[(1.10), (1.14), (1.15), (1.16), (1.17)]{KMM25}, respectively, and
\begin{align*}
	\widetilde{A_{5,n}}(t) &= \widetilde{A_{5,n}}(t,v) \\
	&\coloneqq \frac{i}2 U(t) \int_{1}^{t}  s^{\alpha_0-1}  R_n(s) V_n(s) \\
    & \hspace{20mm} \times \left( \frac{\partial F_n}{\partial z} (v(s))
	\left( i\partial_{t}v(s) - \frac{p^{2}}{2}v(s) \right) - \frac{\partial F_n}{\partial \overline{z}} (v(s))
	\overline{\left( i\partial_{t}v(s) - \frac{p^{2}}{2}v(s) \right)}\right) ds.
\end{align*}
In what follows, we add a tilde to symbols that differ from the definitions in \cite{KMM25}, such as $\widetilde{A_{5,n}}(t)$.

By the standard fixed point argument,
we find a solution to \eqref{E:IE}, instead of \cite[(1.19)]{KMM25}.
Hence the map $\Phi$ as in \cite[(3.1)]{KMM25} is replaced by
\begin{align}\label{E:Phidef}
\begin{aligned}
	\widetilde{\Phi}(v) &\coloneqq U(t-1) v(1) - i {\overline{\lambda_1}} I_1(t,v) \\&\quad -i \sum_{n =2}^\infty  {\overline{\lambda_n}} (A_{1,n}(t) + A_{2,n}(t,v) + A_{3,n}(t,v)  + A_{4,n}(t,v) + \widetilde{A_{5,n}}(t,v)).
\end{aligned}
\end{align}

We further modify the complete metric space $Z_{M}$ as follows:
\begin{align*}
	\widetilde{Z}_{M} &\coloneqq \{ v \in C((0,1],H^1) \mid \norm{v}_X + \norm{v}_{\widetilde{Y}} \le M \}, \qquad M>0, \\
	\norm{v}_{X} &\coloneqq  \norm{ v }_{L^\infty ((0,1];H^{1})} + \norm{v}_{L^{q_0} ((0,1];H^{1}_{r_0})}, \\
	\norm{v}_{\widetilde{Y}} &\coloneqq \norm{\partial_t v}_{L^{q_1} ((0,1];H^{-1}_{r_1})} + \norm{\chi_{\gamma}(x) \mathcal{L} v}_{L^{q_{2}}((0,1]; L^{r_{2}})}, \\
    \chi_{\gamma}(x) &\coloneqq
	\begin{cases}
	|x|^{-\gamma} & \text{if}\ |x| \le 1, \\
	1 & \text{if}\ |x| > 1,
	\end{cases}
    \qquad \mathcal{L} \coloneqq i \partial_{t} - \frac{p^{2}}{2},
\end{align*}
equipped with the metric function
\[
	d_{\widetilde{Z}}(v_1,v_2) \coloneqq \norm{v_1-v_2}_{L^{\infty}((0,1]; L^{2})} + \norm{v_1-v_2}_{L^{q_{0}}((0,1]; L^{r_{0}})}
	+ \norm{\mathcal{L} \left( v_{1} - v_{2} \right)}_{L^{q_{2}}((0,1]; L^{r_{2}})},
\]
where pairs $(q_{j},r_{j})$ are defined as
\begin{align*}
    (q_0,r_0) &\coloneqq \begin{cases}
	(\frac{2}{\alpha_0-1},\frac2{2-\alpha_0}) & d=2, \\
	(4, 3) & d=3,
	\end{cases}
	\qquad
	( q_1 , r_1 ) \coloneqq (1, r_0), \quad  d=2, 3,
    \\
    (q_{2}, r_{2} ) &\coloneqq
    \begin{cases}
    (q_{1}, \rho r_{0}) & d=2, \\
   \left( 1, 2 \right) & d=3,
    \end{cases}
\end{align*}
$3-2 \alpha_0 \le \gamma <1$ and $2 - \alpha_{0} < \rho <1$ if $d=2$, and $\gamma = 1/2$ if $d=3$.
Note that $X$-norm, the pairs $(q_0,r_0)$ and $(q_{1}, r_{1})$ are the same as those in \cite{KMM25}.
To recover from the error, we introduce a truncated weight function $\chi_\gamma$ in the $\widetilde{Y}$-norm,
which behaves like a singular weight $|x|^{-\gamma}$ near the origin, and is truncated to a constant value away from the origin.
This function allows us to localize the singularity near $x = 0$ while maintaining boundedness at $|x| \rightarrow \infty$.
By this truncation, $\chi_{\gamma} \ge 1$ and hence $\norm{\mathcal{L}v}_{L^{q_{2}}(0,1; L^{r_{2}})} \le \norm{\chi_{\gamma}\mathcal{L}v}_{L^{q_{2}}(0,1; L^{r_{2}})}$ holds.
Using the fact and a standard argument (e.g., \cite{K87}), we see that $(\widetilde{Z}_{M}, d_{\widetilde{Z}})$ is a complete metric space for all $M \ge 0$.

The weighted function $\chi_{\gamma}$ is introduced to control the singularity in time $t$.
In two dimensions,
the Duhamel term of \eqref{E:Phiv} exhibits a strong singularity at time $t=0$, which becomes a key issue.
To address this difficulty, we employ resolvent estimates to avoid the time singularity.
However, this remedy leads to spatial singularities at the origin, which makes it necessary to use
the weight function $\chi_{\gamma}$.
In three dimensions, we face a lack of spatial integrability.
To gain spatial integrability, we impose stronger spatial regularity.
This requires us to impose both a time singularity at $t=0$ and a spatial singularity at the origin.
The time singularity is negligible, and the spatial singularity needs to be handled by $\chi_{\gamma}$.
These modifications make it possible to correct the gap in the proof of Theorem 1.4 in \cite{KMM25}.

Here are two remarks for $\widetilde{\Phi}$.
\begin{enumerate}
\item Without using equation \cite[(1.8)]{KMM25}, integration by parts in time in \eqref{E:Phidef} yields
\begin{align}
    \widetilde{\Phi}(v) &= U(t-1) v(1) - i \int_{1}^{t} s^{\alpha_{0}-2} U(t-s) \widetilde{F}(s, v(s))\, ds
    \label{E:Phiv}
\end{align}
for any $v \in \widetilde{Z}_{M}$, where
$\widetilde{F}(t,v) = \sum_{n \ge 1} \overline{ \lambda _ n} e^{-i \frac{(n-1) |x|^2}{2t} } F_n (v)$.
Note that the integration by parts is justified by the argument in \cite[Section 3.3.1]{KMM25}.
We emphasize that
the auxiliary term $\mathcal{E}(v)$, defined by \cite[(3.2)]{KMM25}, is no longer required.
\item For $v \in \widetilde{Z}_{M}$, $\widetilde{\Phi}$ satisfies the equation
\begin{align}\label{E:Phierror}
	i \partial_t \widetilde{\Phi}(v) - \frac{p^2}2 \widetilde{\Phi}(v)
	= t^{\alpha_0-2} \widetilde{F}(t,v).
\end{align}
\end{enumerate}

Under this setting, \cite[Proposition 3.1]{KMM25} is replaced by the following proposition:
\begin{proposition}\label{P:key}
Under the assumption of \cite[Theorem 1.4]{KMM25}, there exist $\varepsilon_0>0$ and $M>0$
such that the following assertions hold:
\begin{enumerate}
\item If $u_0 \in L^2$ satisfies
$u_0 \in J (1) L^2$ and
$\norm{u_0}_{L^2} + \norm{ J(1) u_0 }_{ L^2 } \le \varepsilon_0$ then
$\widetilde{\Phi}$ is a contraction map from $(\widetilde{Z}_M, d_{\widetilde{Z}})$ to itself.
\item If $v \in \widetilde{Z}_{M}$ is a fixed point of
$\widetilde{\Phi}$, then $v$ is a solution to \cite[(1.8)]{KMM25} on $(0,1]$.
\end{enumerate}
\end{proposition}
Thanks to \eqref{E:Phierror}, the second assertion of Proposition \ref{P:key} is immediate.
Before proving the first assertion of Proposition \ref{P:key}, we give the list that
summarizes the changes made to \cite{KMM25}.
Recall that we only need to modify the arguments in two and three dimensions.
\begin{itemize}
    \item Replace $A_{5, n}(t)$ in \cite[(1.18)]{KMM25} by $\widetilde{A_{5,n}}(t)$.
    \item Replace the integral equation \cite[(1.19)]{KMM25} by \eqref{E:IE}.
    \item Replace the map $\Phi$ defined in \cite[(3.1)]{KMM25} by $\widetilde{\Phi}$ defined in \eqref{E:Phidef}.
    \item Remove \cite[(3.2)]{KMM25} and \cite[(3.3)]{KMM25} involving $\mathcal{E}(v)$.
    \item Replace the complete metric space $(Z_{M}, d_{Z})$ defined in \cite[Section 3.1]{KMM25} by $(\widetilde{Z}_{M}, d_{\widetilde{Z}})$.
    \item Replace \cite[Proposition 3.1]{KMM25} by Proposition \ref{P:key}.
    \item Remove \cite[Lemma 3.2]{KMM25} involving $\mathcal{E}(v)$ and its proof.
    \item Replace \cite[Lemma 3.3]{KMM25} by Lemma \ref{L:Phierror} below.
    \item In \cite[Sections 3.2.4 and 3.2.5]{KMM25}, replace the estimates of $A_{5,n}(t)$ by those of $\widetilde{A_{5,n}}(t)$ below.
    \item Replace the proof of \cite[Proposition 3.1]{KMM25} in \cite[Section 3.3]{KMM25} by the proof below.
    \item Replace $Z_{M, T}$ in \cite[Section 3.3.1]{KMM25}
    by $\widetilde{Z}_{M, T}$, which is the same complete metric space as $\widetilde{Z}_{M}$ but time interval is replaced by $(T,1]$.
\end{itemize}

\subsection{Proof of the contraction property}

Let us prove the first assertion of Proposition \ref{P:key}.
Since $\widetilde{\Phi}$ satisfies \eqref{E:Phierror} for any $v \in \widetilde{Z}_{M}$,
the estimate in \cite[Lemma 3.3]{KMM25} is simplified.
In particular, the second term involving $\mathcal{E}(v)$ on the right-hand side is removed as follows:

\begin{lemma} \label{L:Phierror}
Fix $\varepsilon_0>0$ and $M>0$. Then,
\[
	\norm{i \partial_t \widetilde{\Phi}(v) - \frac{p^2}2 \widetilde{\Phi}(v)}_{L^{q_1} ((0,1];H^{-1}_{r_1})}
    \lesssim M^{\alpha}
\]
for $v\in \widetilde{Z}_M$.
\end{lemma}
We shall show
\begin{align}
    \norm{\widetilde{\Phi}(v)}_{X}  + \norm{\widetilde{\Phi}(v)}_{\widetilde{Y}} \le M
    \label{E:Phibo}
\end{align}
for all $v\in \widetilde{Z}_M$.
As for the $X$-norm,
one only treats $\widetilde{A_{5,n}}(t)$, since the other parts are the same as those in \cite{KMM25}.
By Strichartz' estimate, we have
\begin{align}
\begin{aligned}
    \norm{\widetilde{A_{5,n}}(\cdot,v)}_{X} \lesssim{}&
    \norm{s^{\alpha_{0}-1} U(s)
	R_{n}(s) V_{n}(s)
    \frac{\partial F_n}{\partial z} (v(s))
	\mathcal{L}v(s)
    }_{L^{1} ((0,1]; H^{1})} \\
    &+ \norm{ s^{\alpha_{0}-1} U(s)
	R_{n}(s) V_{n}(s)
    \frac{\partial F_n}{\partial \overline{z}} (v(s))
	\overline{\mathcal{L}v(s)} }_{L^{1} ((0,1]; H^{1})}.
\end{aligned}
    \label{E:estA5}
\end{align}
Let us only estimate the first term, because the second term is treated in the same manner.
In what follows, let us continue the proof separately for the two- and three-dimensional cases.

\subsubsection{Two-dimensional case}
When $d = 2$, the singularity in time $t$ in the right-hand side of \eqref{E:estA5} becomes worse.
We hence employ resolvent estimates and instead impose a stronger spatial singularity at the origin to avoid the time singularity.
To deal with the spatial singularity arising at the origin, we introduce the weight function $\chi_\gamma$.
Thanks to $|x|^{- \gamma} \le \chi_{\gamma}$,
we see from \cite[(2.11) and (2.14)]{KMM25} and Gagliardo-Nirenberg's inequality that
\begin{align*}
    &\norm{s^{\alpha_{0}-1} p U(s)
	R_{n}(s) V_{n}(s)
    \frac{\partial F_n}{\partial z} (v(s))
	\mathcal{L}v(s)
    }_{L^{1} ((0,1]; L^{2})} \\
    \lesssim{}& n^{-\frac{1}{2}}
    \norm{s^{\alpha_{0}-\frac{3}{2}+ \frac{\gamma}{2}}
    |x|^{- \gamma} \frac{\partial F_n}{\partial z} (v(s))
	\mathcal{L}v(s)
    }_{L^{1} ((0,1]; L^{2})} \\
    \lesssim{}& n^{\frac{1}{2}} \norm{s^{\alpha_{0}-\frac{3}{2}+ \frac{\gamma}{2}}}_{L^{\infty}((0,1])} \norm{v}_{L^{\infty}((0,1]; L^{\frac{2 \alpha_{0} \rho}{(1- \rho)(2 - \alpha_0)}})}^{\alpha_{0}} \norm{\chi_{\gamma}(x) \mathcal{L} v(s)}_{L^{q_{1}} ((0,1]; L^{\rho r_{0}})} \\
    \lesssim{}& n^{\frac{1}{2}} M^{\alpha}.
\end{align*}
Also, one easily verifies that
\begin{align*}
    \norm{s^{\alpha_{0}-1} U(s)
	R_{n}(s) V_{n}(s)
    \frac{\partial F_n}{\partial z} (v(s))
	\mathcal{L}v(s)
    }_{L^{1} ((0,1]; L^{2})}
    \lesssim n^{\frac{3}{4}} M^{\alpha},
\end{align*}
thereby concluding that
\begin{align}
    \norm{\widetilde{A_{5,n}}(\cdot,v)}_{X} \lesssim n^{\frac{3}{4}} M^{\alpha}.
    \label{E:1}
\end{align}

Let us move on to the $\widetilde{Y}$-norm.
Lemma \ref{L:Phierror} implies that
\begin{align*}
    \norm{\partial_t \widetilde{\Phi}(v)}_{L^{q_1} ((0,1];H^{-1}_{r_1})}
    &\le \frac{1}{2} \norm{\widetilde{\Phi}(v)}_{L^{q_0} ((0,1];H^{1}_{r_1})}
    + \norm{i \partial_t \widetilde{\Phi}(v) - \frac{p^2}2 \widetilde{\Phi}(v)}_{L^{q_1} ((0,1];H^{-1}_{r_1})} \\
    &\lesssim \norm{\widetilde{\Phi}(v)}_{X} + M^{\alpha}
\end{align*}
for all $v \in \widetilde{Z}_{M}$.
Furthermore, by \cite[Lemma 2.12]{KMM25} and Gagliardo-Nirenberg's inequality, we see from \eqref{E:Phierror} that
\begin{align*}
    &\norm{\chi_{\gamma}(x) \mathcal{L} \widetilde{\Phi}(v)}_{L^{q_{1}} ((0,1]; L^{\rho r_{0}})} \\
    \le{}& \sum_{n \ge 1} \abs{\lambda_{n}} \norm{t^{\alpha_{0}-2} |v|^{\alpha_{0}} \chi_{\gamma}(x) v}_{L^{q_{1}} ((0,1]; L^{\rho r_{0}})} \\
    \lesssim{}& \norm{t^{\alpha_{0}-2}}_{L^{\frac{2}{3- \alpha_{0}}}((0,1])} \norm{v}_{L^{\infty} ((0,1]; L^{\frac{2 \alpha_{0} \rho}{(1- \rho)(2 - \alpha_0)}})}^{\alpha_{0}} \\
    &\quad
    \times \left( \norm{|x|^{-\gamma} v}_{L^{q_{0}}((0,1]; L^{r_{0}}(|x|\le 1))} + \norm{v}_{L^{q_{0}}((0,1]; L^{r_{0}}(|x| > 1))} \right) \\
    \lesssim{}& \norm{v}_{L^{\infty}((0,1]; H^{1})}^{\alpha_{0}} \norm{v}_{L^{q_{0}}((0,1]; H^{1}_{r_{0}})}
    \lesssim M^{\alpha}
\end{align*}
for all $v \in \widetilde{Z}_{M}$.
Combining the above estimates, we obtain \eqref{E:Phibo}
for all $v\in \widetilde{Z}_M$ if $M$ is small enough.

Let us prove the contraction property
\begin{align}
    d_{\widetilde{Z}} (\widetilde{\Phi}(v_1),\widetilde{\Phi}(v_2)) \le \frac12 d_{\widetilde{Z}} (v_1,v_2)
    \label{E:con1}
\end{align}
for all $v_1,v_2 \in \widetilde{Z}_M$ if $M$ is small enough.
Thanks to \eqref{E:Phiv} and \eqref{E:Phierror}, the contraction property can be shown by the standard well-posedness theory.
Indeed, using
\begin{align}
    \abs{\widetilde{F}(t, v_1) - \widetilde{F}(t, v_2)} \lesssim \left(|v_1|^{\alpha-1} + |v_2|^{\alpha-1}  \right) |v_1 - v_2|,
    \label{non:e1}
\end{align}
and Hardy's inequality, similarly to the above, we obtain
\begin{align*}
    &\norm{\mathcal{L} \left( \widetilde{\Phi}(v_{1}) - \widetilde{\Phi}(v_{2}) \right)}_{L^{q_{2}}((0,1]; L^{r_{2}})} \\
    \lesssim{}& \norm{t^{\alpha_{0}-2} \left( \widetilde{F}(t, v_1) - \widetilde{F}(t, v_2) \right)}_{L^{1}((0,1]; L^{r_{0}})} \\
    \lesssim{}& \norm{t^{\alpha_{0}-2}}_{L^{\frac{2}{3- \alpha_{0}}}((0,1])}
    \left( \norm{v_{1}}_{L^{\infty} ((0,1]; L^{\frac{2 \alpha_{0} \rho}{(1- \rho)(2 - \alpha_0)}})}^{\alpha_{0}} + \norm{v_{2}}_{L^{\infty} ((0,1]; L^{\frac{2 \alpha_{0} \rho}{(1- \rho)(2 - \alpha_0)}})}^{\alpha_{0}} \right)
    \norm{v_{1} - v_{2}}_{L^{q_{0}}((0,1]; L^{r_{0}})} \\
    \lesssim{}& M^{\alpha_{0}} \norm{v_1 - v_2}_{L^{q_{0}}((0,1]; L^{r_{0}})}
\end{align*}
for any $v_1$, $v_2 \in \widetilde{Z}_{M}$.
Note that since the metric $d_{\widetilde{Z}}$ does not involve the truncated weight $\chi_{\gamma}$, the above estimate can be closed in this metric.
Further, by Strichartz' estimate, one has
\begin{align*}
    &\norm{\widetilde{\Phi}(v_1) - \widetilde{\Phi}( v_2)}_{L^{\infty}((0,1]; L^{2})} + \norm{\widetilde{\Phi}(v_1) - \widetilde{\Phi}( v_2)}_{L^{q_{0}}((0,1]; L^{r_{0}})} \\
    \lesssim{}& \norm{ s^{\alpha_0 -2} ( |v_1(s)|^{\alpha_{0}}+ |v_{2} (s)|^{\alpha_{0}} ) }_{L^{q_{3}}((0,1] ; L^{r_{3}})}
    \left( \norm{v_{1} - v_{2}}_{L^{\infty}((0,1]; L^{2})} + \norm{v_{1} - v_{2}}_{L^{q_{0}}((0,1]; L^{r_{0}})}  \right),
\end{align*}
where $(q_{3}, r_{3}) = \left( 2/(\alpha_{0} -1), 2/(3- \alpha_{0}) \right)$.
A use of Gagliardo-Nirenberg inequality leads to
\begin{align*}
    &\norm{ s^{\alpha_{0} -2} ( |v_1(s)|^{\alpha_{0}}+ |v_{2} (s)|^{\alpha_{0}} ) }_{L^{q_{3}}((0,1] ; L^{r_{3}})} \\
    \le{}& \norm{ s^{\alpha_{0} -2}}_{L^{q_{3}}((0,1])}
    \left( \norm{v_1(s)}_{L^{\infty}((0,1] ; L^{\alpha_{0}r_{3}})}^{\alpha_{0}} + \norm{v_2(s)}_{L^{\infty}((0,1] ; L^{\alpha_{0}r_{3}})}^{\alpha_{0}} \right)
    \lesssim M^{\alpha-1}
\end{align*}
for any $v_1$, $v_2 \in \widetilde{Z}_{M}$, which implies that \eqref{E:con1} holds if $M$ is small enough.

\subsubsection{Three-dimensional case}

We shall handle the case $d=3$.
Recall that $\gamma =1/2$.
In this case, we face a lack of spatial integrability.
To address this issue, we employ Sobolev embedding
\[
    \norm{f}_{L^{2}} \lesssim \norm{p^{\frac{1}{2}}f}_{L^{\frac{3}{2}}}.
\]
One also utilizes the resolvent estimate to handle the additional spatial regularity. This requires imposing an additional singularity at time $t=0$ and at the spatial origin.
While the time singularity is negligible,
the spatial singularity needs to be handled by the weight function $\chi_\gamma$ as follows:
Together with Sobolev embedding and $|x|^{-1/2} \le \chi_{\gamma}$,
we see from \cite[(2.11) and (2.14)]{KMM25} that
\begin{align}
    \begin{aligned}
    \norm{p U(s) R_{n}(s) V_{n}(s) f g}_{L^{2}}
    &\lesssim \norm{p^{\frac{3}{2}} U(s) R_{n}(s) V_{n}(s) f g}_{L^{3/2}} \\
    &\lesssim s^{-\frac{3}{4}} \norm{|\sqrt{s} p|^{\frac{3}{2}} U(s) R_{n}(s) V_{n}(s) \abs{\frac{x}{\sqrt{s}}}^{\frac{1}{2}} \abs{\frac{x}{\sqrt{s}}}^{- \frac{1}{2}} f g}_{L^{3/2}} \\
    &\lesssim s^{- \frac{1}{2}} n^{-\frac{1}{4}} \norm{f}_{L^{6}} \norm{\chi_{1/2}(x) g}_{L^{2}}
    \end{aligned}
    \label{E:3d}
\end{align}
for any $f \in L^{6}$ and $g \in L^{2}$.
Recalling $\alpha_{0} = 3/2$, Gagliardo-Nirenberg's inequality and \eqref{E:3d} lead to
\begin{align*}
    &\norm{s^{\alpha_{0}-1} p U(s)
	R_{n}(s) V_{n}(s)
    \frac{\partial F_n}{\partial z} (v)
	\mathcal{L}v
    }_{L^{1} ((0,1]; L^{2})} \\
    \lesssim{}& n^{\frac{3}{4}}
    \norm{v}_{L^{\infty} ((0,1]; L^{6})}
	\norm{\chi_{1/2}(x) \mathcal{L}v
    }_{L^{1} ((0,1]; L^{2})} \\
    \lesssim{}& n^{\frac{3}{4}}
    \norm{v}_{L^{\infty} ((0,1]; H^{1})}
	\norm{\chi_{1/2}(x) \mathcal{L}v
    }_{L^{1} ((0,1]; L^{2})} \\
    \lesssim{}& n^{\frac{3}{4}} M^{\alpha}.
\end{align*}
Also, we easily obtain
\begin{align*}
    &\norm{s^{\alpha_{0}-1} U(s)
	R_{n}(s) V_{n}(s)
    \frac{\partial F_n}{\partial z} (v)
	\mathcal{L}v
    }_{L^{1} ((0,1]; L^{2})} \\
    \lesssim{}& n^{\frac{7}{8}} \norm{s^{\frac{1}{4}}}_{L^{\infty}(0,1)}
    \norm{v}_{L^{\infty} ((0,1]; L^{6})}
	\norm{\mathcal{L}v
    }_{L^{1} ((0,1]; L^{2})} \\
    \lesssim{}& n^{\frac{7}{8}} M^{\alpha}.
\end{align*}

Hence, we conclude that
\begin{align}
    \norm{\widetilde{A_{5,n}}(\cdot,v)}_{X} \lesssim n M^{\alpha}.
    \label{E:1-3d}
\end{align}

Let us treat the $\widetilde{Y}$-norm.
The estimate
\begin{align*}
    \norm{\partial_t \widetilde{\Phi}(v)}_{L^{q_1} ((0,1];H^{-1}_{r_1})}
    &\lesssim \norm{\widetilde{\Phi}(v)}_{X} + M^{\alpha}
\end{align*}
is the same as in the two-dimensional case; we therefore focus on the second term.
By \cite[Lemma 2.12]{KMM25} and Gagliardo-Nirenberg's inequality, one sees from \eqref{E:Phierror} that
\begin{align*}
    &\norm{\chi_{\gamma}(x) \mathcal{L} \widetilde{\Phi}(v)}_{L^{q_{2}} ((0,1]; L^{r_{2}})} \\
    \le{}& \sum_{n \ge 1} \abs{\lambda_{n}} \norm{t^{\alpha_{0}-2} |v| \chi_{1/2}(x) v}_{L^{1} ((0,1]; L^{2})} \\
    \lesssim{}& \norm{t^{-\frac{1}{2}}}_{L^{\frac{4}{3}}((0,1])} \norm{v}_{L^{\infty} ((0,1]; L^{6})} \left(
    \norm{|x|^{-\frac{1}{2}} v}_{L^{q_{0}}((0,1]; L^{r_{0}}(|x| \le 1))}
    +
    \norm{v}_{L^{q_{0}}((0,1]; L^{r_{0}}(|x|>1))}
    \right) \\
    \lesssim{}& \norm{v}_{L^{\infty}((0,1]; H^{1})} \norm{v}_{L^{q_{0}}((0,1]; H^{1}_{r_{0}})}
    \lesssim M^{\alpha}
\end{align*}
for all $v \in \widetilde{Z}_{M}$.
Combining the above estimates, we obtain \eqref{E:Phibo}
for all $v\in \widetilde{Z}_M$ if $M$ is small enough.

Let us prove the contraction property \eqref{E:con1} for all $v_1,v_2 \in \widetilde{Z}_M$ if $M$ is small enough.
Arguing as in the two-dimensional case, using \eqref{non:e1}, we have
\begin{align*}
    &\norm{\mathcal{L} \left( \widetilde{\Phi}(v_{1}) - \widetilde{\Phi}(v_{2}) \right)}_{L^{q_{2}}((0,1]; L^{r_{2}})} \\
    \lesssim{}& \norm{t^{-\frac{1}{2}}}_{L^{ \frac{4}{3}}((0,1])}
    \left( \norm{v_{1}}_{L^{\infty} ((0,1]; L^{6})} + \norm{v_{2}}_{L^{\infty} ((0,1]; L^{6})} \right)
    \norm{v_{1} - v_{2}}_{L^{q_{0}}((0,1]; L^{r_{0}})} \\
    \lesssim{}& M \norm{v_1 - v_2}_{L^{q_{0}}((0,1]; L^{r_{0}})}
\end{align*}
for any $v_1$, $v_2 \in \widetilde{Z}_{M}$.
Strichartz' estimate and Gagliardo-Nirenberg's inequality lead to
\begin{align*}
    &\norm{\widetilde{\Phi}(v_1) - \widetilde{\Phi}( v_2)}_{L^{\infty}((0,1]; L^{2})} + \norm{\widetilde{\Phi}(v_1) - \widetilde{\Phi}( v_2)}_{L^{q_{0}}((0,1]; L^{r_{0}})} \\
    \lesssim{}& \norm{ s^{-\frac{1}{2}} ( |v_{1}| + |v_{2}| ) }_{L^{4/3}((0,1] ; L^{6})}
    \left( \norm{v_{1} - v_{2}}_{L^{\infty}((0,1]; L^{2})} + \norm{v_{1} - v_{2}}_{L^{q_{0}}((0,1]; L^{r_{0}})}  \right) \\
    \lesssim{}& \norm{ s^{-\frac{1}{2}}}_{L^{4/3}((0,1])}
    \left( \norm{v_{1}}_{L^{\infty}((0,1] ; L^{6})} + \norm{v_{2}}_{L^{\infty}((0,1] ; L^{6})} \right)
    d_{\widetilde{Z}} (v_1,v_2)
    \lesssim M d_{\widetilde{Z}} (v_1,v_2)
\end{align*}
for any $v_1$, $v_2 \in \widetilde{Z}_{M}$, which implies that \eqref{E:con1} holds if $M$ is small enough.

\bigskip

Based on the preceding subsections,
we conclude that the map $\widetilde{\Phi}$ is a contraction on $(\widetilde{Z}_{M}, d_{\widetilde{Z}})$ for sufficiently small $M > 0$.
This completes the proof of the first assertion of Proposition \ref{P:key}.
\qed

\section*{Acknowledgments}
M.K. was supported by JSPS KAKENHI Grant Number 24K06796.
S.M. was supported by JSPS KAKENHI Grant Number 24K00529.
H.M. was supported by JSPS KAKENHI Grant Number 22K13941.


\end{document}